\documentclass[review]{elsarticle}
\usepackage{amsfonts}
\usepackage{amsmath}
\usepackage{amssymb}
\usepackage{amsthm}   
\usepackage{empheq}   
\usepackage{titlesec} 
\allowdisplaybreaks

\newtheorem{thm}{Theorem}[section]
\newtheorem{lem}{Lemma}[section]
\newtheorem{rem}{Remark}[section]
\usepackage{bm}
\usepackage{bookmark}
\usepackage{arydshln}
\usepackage{lineno,hyperref}
\modulolinenumbers[5]

\journal{Journal of \LaTeX\ Templates}
\begin{document}
\begin{frontmatter}

\title{Long-time existence of nonlinear inhomogeneous compressible elastic waves}

\author[a]{Silu Yin\corref{mycorrespondingauthor}}
\cortext[mycorrespondingauthor]{Corresponding author}
\ead{yins11@shu.edu.cn}

\author[b]{Xiufang Cui}
\ead{xfcui16@fudan.edu.cn}

\address[a]{Department of Mathematics, Shanghai University, Shanghai, 200444, P. R. China}
\address[b]{School of Mathematical Sciences, Fudan University, Shanghai, 200433, P. R. China}

\begin{abstract}
 In this paper, we consider the nonlinear inhomogeneous compressible elastic waves in three spatial dimensions when the density is a small disturbance around a constant state. In homogeneous case, the almost global existence was established by Klainerman-Sideris \cite{KS}, and global existence was built by Agemi \cite{Ag} and Sideris \cite{S1,S2} independently. Here we establish the corresponding almost global and global existence theory in the inhomogeneous case.
\end{abstract}


\end{frontmatter}

\linenumbers
\setcounter{section}{0}
\numberwithin{equation}{section}

\section{Introduction}\label{s1}
The motion of an elastic body in three spatial dimensions is described by a time-dependent family of orientation preserving diffeomorphisms, written as $\varphi=\varphi(t,x),\ 0\leq t<T$, where $\varphi(0,x)=x$. Here $x=(x^1,x^2,x^3)$ is called material point, $\varphi=(\varphi^1,\varphi^2,\varphi^3)$ is called spatial point. The deformation gradient tensor $F$ is denoted as
  $$F=\frac{\partial \varphi(t, x)}{\partial x}\Big|_{x=x(t,\varphi)},$$
  where $F^j_i=\frac{\partial \varphi^j}{\partial x^i}$.

We concentrate on three dimensional inhomogeneous compressible hyperelastic materials whose potential energy $W(F)$ is determined only by the deformation tensor $F$. Denote the initial density by $\rho(x)$, then the Lagrangian variation governing the  elastodynamics materials takes
\begin{align}\label{1.1}
  \mathcal{\delta}(\varphi)=\int\int_{\mathbb{R}^3}\Big(\frac12\rho(x)|\partial_t\varphi(t,x)|^2-W(F)\Big)dxdt.
\end{align}
Then, we can formulate \eqref{1.1} to a nonlinear system
\begin{align}\label{11}
  \rho(x)\partial_t^2\varphi-\nabla_x\cdot\frac{\partial W(\nabla_x \varphi)}{\partial F}=0.
\end{align}

From the groundbreaking work of quasilinear wave equations and elastodynamics in three dimensions, both the smallness of initial data and null condition in nonlinearities are indispensable to make sure the solution exists globally. It is reasonable to consider the elastodynamics around the equilibrium state. Let $u(t,x)=\varphi(t,x)-x$, by the isotropy of the hyperelastic materials and inherent frame-indifference, the mechanism obeys
\begin{align}\label{1.2}
 \rho(x)\partial_t^2 u-c_2^2\Delta u-(c_1^2-c_2^2)\nabla( \nabla\cdot u )= N(u,u)
\end{align}
deduced from \eqref{11} , where $c_1^2>\frac43c_2^2$, and
$$N^i(u,v)=B_{lmn}^{ijk}\partial_l(\partial_mu^j\partial_nu^k)$$
 with symmetry
 \begin{align}\label{1.4}
B_{lmn}^{ijk}=B_{mln}^{jik}=B_{lnm}^{ikj}.
 \end{align}
Here we use the null condition appeared in \cite{S2}, which means that the self-interaction of each homogeneous elastic wave family is nonresonant:
\begin{equation}\label{1.5}
B_{lmn}^{ijk}\omega_i\omega_j\omega_k\omega_l\omega_m\omega_n=0,\quad for \ all\  \omega\in S^2
\end{equation}
\begin{equation}\label{1.6}
B_{lmn}^{ijk}\eta_i\eta_j\eta_k\omega_l\omega_m\omega_n=0, \quad  for\  all \  \eta, \omega \in S^2 \ \textup{with} \  \eta\bot\omega.
\end{equation}

For three dimensional homogeneous compressible elastic waves in which the density is always a constant, John \cite{J2} first showed a small displacement solution exists almost globally via an $L^1-L^\infty$ estimate (see also \cite{KS} dependently by the generalized energy method). Sideris \cite{S1} discovered a null condition within the class of physically meaningful nonlinearities and established the global existence theory via an enhanced decay estimate. Refinements for the proof of global existence were presented by Sideris \cite{S2} and Agemi \cite{Ag} independently.  Blow up for elasticity with large initial data was given by Tahvildar-Zadeh \cite{T}. The global well-posedness of compressible elastic waves in two dimensions is more delicate and still remains open.

In this paper, we consider the long-time existence of inhomogeneous compressible elastic waves in which the density is not a constant any more. By adapting the improved generalized vector theory of Klainerman \cite{K2} and inspired by Sideris \cite{S2}, we show a unique solution of \eqref {1.2} exists almost globally for non-degenerate elastodynamics near the equilibrium state; see Theorem \ref{thm 2.1}. We also prove that if the elastodynamics is degenerate in the nonlinearities, that is the system satisfies the null condition as in \eqref{1.5} and \eqref{1.6}, then the solution exists globally; see Theorem \ref{thm 2.2}. For inhomogeneous case, not only the Lorentz invariance is absent but also there will be some extra linear terms to be controlled. One difference and difficulty we meet here is that the traditional weighted norm $\mathcal{X}_{|\alpha|+2}$ does not include $\|\langle c_at-r\rangle P_a\partial_t^2 \Gamma^\alpha u\|_{L^2}$. Utilizing the equation, we discover that $\|\langle c_at-r\rangle P_a\partial_t^2 \Gamma^\alpha u\|_{L^2}$ can also be controlled by the classical generalized energy $\mathcal{E}_{|\alpha|+2}$; see Lemma \ref{lem43}.

 Next we highlight some closely related results. For homogeneous incompressible elastodynamics, the equations are inherently linearly degenerate in the isotropic case and satisfy a null condition that is necessary for global existence in three dimensions (see \cite{STh1, STh2}). In those papers the authors showed the lower-order energy is bounded but the higher-order energy has a low-increase ratio in time. Recently, Lei-Wang \cite{LW} proved that the higher-order generalized energy is still uniformly bounded and  provided an improved proof for the global well-posedness.
The difficulty in two dimensional case is highly essential. The first nontrivial long-time existence result concerning on this problem was established by Lei-Sideris-Zhou \cite{LSZ}. They proved the almost global existence by formulating the system in Euler coordinate. The breakthrough of global well-posedness was owing to Lei \cite{L}. He introduced the so-called notion of "strong null condition" and discovered the two dimensional incompressible elastodynamics inherently satisfies this extremely important structure in Lagrangian coordinate. Based on these discoveries, Lei proved the incompressible isotropic elastodynamics in two dimensions admits a unique global classical solution. The method can even be used to prove the vanishing viscosity limit of viscoelastic system, see \cite{Ca}. Wang \cite{W} built the global existence theory of incompressible elastodynamics in two dimensions in frequency space in Euler coordinate. For the inhomogeneous case, Yin \cite{yin} built a model of two-dimensional incompressible elastodynamics and showed the global existence. The global well-posedness theory of inhomogeneous incompressible case in three dimensions is coming out in our another paper.

For three dimensional wave equations, John \cite{JF} first proved the blow up phenomenon for the Cauchy problem of wave equations with sufficiently small initial data but violating null condition. John-Klainerman \cite{JK} showed the almost global existence theory for nonlinear scalar wave equations. Then Klainerman \cite{K2} proved the global existence of classical solutions. This crucial work was also obtained independently by Christodoulou \cite{C} using a conformal mapping method. The multiple-speeds case in three dimensions was achieved by Sideris-Tu \cite{ST} under a null condition. In two-dimensional scalar case, Alinhac established a series of results, a "blow-up solution of cusp type" with sufficiently small initial data but without null condition in \cite{Alin1, Alin2} and the global existence with null bilinear forms in \cite{Alin3, Alin4} were obtained. However, the well-posedness of nonlinear wave systems with multiple speeds in two dimensions is still unknown.

Before ending this section, we make an outline of this paper. We first make some preparatory work and state our main results in Section 2. In Section 3, the detailed commutation relationships for the inhomogeneous compressible elastodynamics will be considered. We present some useful $L^\infty-L^2$ estimates and weighted $L^2-L^2$ estimates based on the classical Sobolev inequalities without Lorentz operators in Section 4.  We control the weighted generalized $L^2$ norm of $\partial_t^2 u$ by the classical generalized energy. In Section 5, the almost global existence of inhomogeneous compressible elastodynamics will be established by a higher-order weighted generalized energy estimate. The relative lower-order energy estimate will be presented in Section 6 under the null condition. Combined with the higher-order energy estimate, we can obtain the global existence theory of the inhomogeneous compressible elastodynamics.

\section{Notations and Main Results}

In this work, we concentrate on the long-time existence of the solutions of inhomogeneous elastodynamics where the density is a small disturbance around a constant state. Without loss of generality, let this constant state be $1$. Thus, we assume
$$\rho(x)=1+\tilde{\rho}(x).$$

Partial derivatives will be presented as
\begin{align*}
\partial = (\partial_0,\partial_1,\partial_2,\partial_3)=(\partial_t,\partial_1,\partial_2,\partial_3),\quad \nabla=(\partial_1,\partial_2,\partial_3).
\end{align*}
The angular momentum operators are the vector fields
$$ \Omega=(\Omega_1,\Omega_2,\Omega_3)=x\wedge\nabla, $$
$$\nabla=\frac{x}{r}\partial_r-\frac{x}{r^2}\wedge \Omega ,  \quad \textup{where} \quad r=|x| \quad \textup{and}  \quad \partial_r=\frac{x}{r}\cdot\nabla .$$
By the generators of simultaneous rotations and scaling transform which can be referred to Sideris \cite{S1,S2}, the angular momentum operators  and scaling operator are given as
 \begin{align*}
 \tilde{\Omega}_l=\Omega_lI+U_l,\quad l=1,2,3,
 \end{align*}
 with
\begin{equation*}
U_1=\left[
\begin{matrix}
 0 &0 & 0\\
 0 & 0 &  1 \\
0&-1&0
\end{matrix}
\right],
\quad
U_2=\left[
\begin{matrix}
 0 &0 & -1\\
 0 & 0 &  0 \\
1 &0 &0
\end{matrix}
\right],
\quad
U_3=\left[
\begin{matrix}
 0 &1 & 0 \\
 -1 & 0 &  0 \\
0 &0 &0
\end{matrix}
\right],
\end{equation*}
and
\begin{align*}
 \tilde{S}=t\partial_t +r\partial_r-1.
\end{align*}

Let
\begin{align*}
  \Gamma = (\Gamma_0,\cdots,\Gamma_7)=(\partial,\tilde{\Omega},\tilde{S}),\quad \textup{where} \quad \tilde{\Omega}=\{\tilde{\Omega}_1,\tilde{\Omega}_2,\tilde{\Omega}_3\}.
\end{align*}

As in \cite{S2}, the standard generalized energy is defined as
\begin{equation*}
  \mathcal{E}_k(u(t)) = \frac{1}{2} \sum\limits_{|\alpha|\leq k-1}\int_{\mathbb{R}^3}[ \ |\partial_t \Gamma ^ \alpha u(t)|^2 + c_2^2|\nabla \Gamma ^ \alpha u(t)|^2+(c_1^2-c_2^2)(\nabla \cdot \Gamma ^ \alpha u(t))^2 \ ]dx.
\end{equation*}
To estimate the different family of elastic waves, orthogonal projections onto radial and transverse directions are introduced as follows
\begin{align*}
P_1u(x)=\frac{x}{r}\otimes\frac{x}{r} u(x)=\frac{x}{r}\langle\frac{x}{r},u(x)\rangle\quad
\textup{and}\quad
P_2u(x)=-\frac{x}{r}\wedge(\frac{x}{r}\wedge u(x)).
\end{align*}
And the weighted generalized energy is given by
\begin{equation*}
\mathcal{X}_k(u(t)) =\sum_{a =1}^{2}\sum_{\beta=0}^{3}\sum_{l=1}^{3}\sum_{|\alpha|\leq k-2} ||\langle c_a t-r \rangle P_a \partial_\beta \partial_l\Gamma^\alpha u(t)||_{L^2}
\end{equation*}
with the notation of $\langle \cdot \rangle = (1+|\cdot|^2)^{\frac{1}{2}}$.
We characterize the space of initial data in $H_\Lambda^k$, which is defined as
$$H_\Lambda^k=\{(f,g):\sum\limits_{|\alpha|\leq k-1}(\|\Lambda^\alpha f\|_{L^2}+\|\nabla\Lambda^\alpha f\|_{L^2}+\|\Lambda^\alpha g\|_{L^2})<\infty\},$$
where $$\Lambda=\{\triangledown,\tilde{\Omega},r\partial_r-1\}.$$
For simplicity of presentation, throughout this paper, we utilize $A \lesssim B $ to denote $A\leq CB$ for some positive absolute constant $C$.

Now we are ready to state our main results, which generalize the work of Klainerman-Sideris \cite{KS}, Agemi \cite{Ag} and Sideris \cite{S1,S2} to the inhomogeneous case.
\begin{thm}[Almost Global Existence]\label{thm 2.1}
 Let $k\geq 9$. Suppose $U_0=(u(0),u_t(0))$ is the initial data of \eqref{1.2} and satisfies
 \begin{align*}
 ||U_0||_{H_\Lambda^k}\leq M, \quad  \textup {and}\quad  ||U_0||_{H_\Lambda^{k-2}}\leq \epsilon,
 \end{align*}
where $M,\epsilon >0$ are two given constants.
If
\begin{equation}\label{21}
\|\langle r\rangle\Lambda^\alpha\tilde{\rho}(x)\|_{L^2}\leq\delta<\frac12,\quad\textup{for}\ |\alpha|\leq k+2
\end{equation}
is small enough, then there exists a $\epsilon_0$ sufficient small, which only depends on $M, k, \delta$ such that for any $\epsilon \leq \epsilon_0$, the Cauchy problem of \eqref{1.2} has a unique almost global solution satisfying
\begin{equation}
  \mathcal{E}_{k}^\frac12(u(t))\leq C M\langle t \rangle^{(\varepsilon+\delta)/2}
\end{equation}
for some positive constant $C$ uniformly in $t$.
\end{thm}

\begin{thm}[Global Existence]\label{thm 2.2}
Under the assumption of Theorem \ref{thm 2.1}, if \eqref{1.5} and \eqref{1.6} are satisfied, then the Cauchy problem of  \eqref{1.2} admits a unique global classical solution satisfying
\begin{equation}
  \mathcal{E}_{k-2}^\frac12(u(t))\leq C_1(\varepsilon\exp\{C_2M\}+\delta \exp\{C_3M\})
\end{equation}
for some positive constant $C_1$, $C_2$, $C_3$ uniformly in $t$.
\end{thm}

\begin{rem}
In Theorem \ref{thm 2.1} and Theorem \ref{thm 2.2}, we need some decay in $\langle r\rangle$ for the disturbance $\tilde{\rho}$ of density. In fact, this is not a necessary request. If the density is a small disturbance with a compact support, Theorem \ref{thm 2.1} and Theorem \ref{thm 2.2} are also valid. By searching some weaker condition about the density to keep the long-time existence of the inhomogeneous elastic waves, it is easy to see that assumption \eqref{21} can be replaced by
$$\|\langle r\rangle\Lambda^\alpha\tilde{\rho}(x)\|_{L^2}\leq\delta<\frac12,\quad\textup{for}\ |\alpha|\leq 4$$
and
$$\|\Lambda^\alpha\tilde{\rho}(x)\|_{L^2}\leq\delta<\frac12,\quad\textup{for}\ 5\leq|\alpha|\leq k$$
from the process of proofs.
\end{rem}
\section{Commutation}
To do the generalized energy estimate, it is necessary to analyze what happens if $\Gamma$ derivatives act on \eqref{1.2}. We formulate \eqref{1.2} to
\begin{align}\label{31}
\partial_t^2 u-c_2^2\Delta u -(c_1^2-c_2^2)\nabla (\nabla \cdot u)=N(u,u)-\tilde{\rho}\partial_t^2 u.
\end{align}
Define
\begin{equation}
\mathcal{L}u\triangleq\partial_t^2 u-\mathcal{A}u\triangleq\partial_t^2 u -c_2^2\Delta u -(c_1^2-c_2^2)\nabla(\nabla \cdot u).
\end{equation}
By direct calculation, we obtain the following commutation
\begin{align*}
  \partial\mathcal{L}u=\mathcal{L}\partial u,\quad \tilde{\Omega}\mathcal{L}u=\mathcal{L}\tilde{\Omega} u,\quad \tilde{S}\mathcal{L}u=\mathcal{L}\tilde{S} u-2\mathcal{L}u,\quad \tilde{S}\partial_t^2u=\partial_t^2\tilde{S} u-2\partial_t^2u,
\end{align*}
and
\begin{align*}
\partial N(u,v)=&N(\partial u,v)+N( u,\partial v),& \partial(\tilde{\rho}\partial_t^2 u)&=\nabla\tilde{\rho}\partial_t^2 u+\tilde{\rho}\partial_t^2 \partial u,\\
\tilde{\Omega}N(u,v)=&N(\tilde{\Omega} u,v)+ N(\tilde{\Omega} u ,v),& \tilde{\Omega}(\tilde{\rho}\partial_t^2 u)&=\tilde{\Omega}\tilde{\rho}\partial_t^2 u+\tilde{\rho}\partial_t^2 \tilde{\Omega} u,\\
\tilde{S}N(u,v)=&N(\tilde{S} u,v)+N( u,\tilde{S}v)-2N(u,v),&\tilde{S}(\tilde{\rho}\partial_t^2 u)&=(r\partial_r-1)\tilde{\rho}\partial_t^2 u+\tilde{\rho}\partial_t^2 \tilde{S} u-2\tilde{\rho}\partial_t^2 u.
\end{align*}
Then we have
\begin{align*}
\mathcal{L}\partial u=&2N(\partial u,u)-\nabla\tilde{\rho}\partial_t^2 u-\tilde{\rho}\partial_t^2 \partial u,\\
\mathcal{L}\tilde{\Omega} u=&2N(\tilde{\Omega} u,u)-\tilde{\Omega}\tilde{\rho}\partial_t^2 u-\tilde{\rho}\partial_t^2 \tilde{\Omega} u,\\
\mathcal{L}\tilde{S}u=&2N(\tilde{S} u,u)-(r\partial_r-1)\tilde{\rho}\partial_t^2 u+\tilde{\rho}\partial_t^2 \tilde{S} u
\end{align*}
according to symmetry \eqref{1.4}. For any multi-index $\alpha=(\alpha_1,\cdots,\alpha_8)\in\mathbb{N}^8$, we then can get the following
\begin{equation}\label{3.1}
\mathcal{L}\Gamma^\alpha u
\backsimeq\sum_{|\beta+\gamma|= |\alpha|}N(\Gamma^\beta u,\Gamma^\gamma u)-\sum_{|\beta+\gamma|= |\alpha|}\Lambda^\beta\tilde{\rho}\partial_t^2\Gamma^\gamma u.
\end{equation}
Here we use notation $A\backsimeq E+F$ means that $A=k_1 E+k_2F$ for some finite positive numbers $k_1$ and $k_2$ which depend only on $\alpha$.

\section{$L^\infty$ and weighted $L^2$ estimates}
From the theory of nonlinear wave equations, we know that a priori estimates play the key role of global existence. Klainerman \cite{K} obtained pointwise bounds for the unknown that decay as $t\rightarrow\infty$ by building some Klainerman-Sobolev estimates with a larger collection of vector fields that preserve the linear wave equation. Armed with these estimates, it is not hard to adapt the proof of the local existence theorem to obtain global existence in certain cases.  But for elastodynamics, the Lorentz invariance is absent. The following improved Sobolev-type estimates were appeared partially in  \cite{S2} (see also \cite{KS}).
\begin{lem}[\cite{S2}]\label{lem4.1}
Let $u\in H_\Gamma^k(T)$  and  $\mathcal{X}_k (u(t))<\infty, $ then
\begin{align}
\langle r \rangle^{\frac{1}{2}}| \Gamma^\alpha u(t,x)|&\leq C\mathcal{E}_k^\frac{1}{2}(u(t)), \quad|\alpha|+2 \leq k,\label{4.1}\\
\langle r \rangle|\partial \Gamma^\alpha u(t,x)|&\leq C\mathcal{E}_k^\frac{1}{2}(u(t)),  \quad|\alpha|+3 \leq k,\label{4.2-1}\\
\langle r \rangle {\langle c_\alpha t -r\rangle}^{\frac{1}{2}}|P_\alpha \partial\Gamma^\alpha u(t,x)|&\leq C[\mathcal{E}_k^{\frac{1}{2}}(u(t))+\mathcal{X}_k (u(t))],\quad |\alpha|+3\leq k,\label{4.2}\\
 \langle r \rangle\langle c_\alpha t -r \rangle |P_\alpha \partial \nabla \Gamma ^\alpha u(t,x)|&\leq C \mathcal{X}_k(u(t)), \quad |\alpha |+4 \leq k.\label{4.3}
\end{align}
\end{lem}
For small solutions of the inhomogeneous elastodynamics, the weighted norm $\mathcal{X}_k$ can be controlled by the energy $\mathcal{E}_k^\frac12$. Before doing this, we first state a basic estimate appeared in \cite{S2}.
\begin{lem}\label{lem41}
  Let $ u\in H_\Gamma^2(T) $, then
 \begin{align*}
  \mathcal{X}_2(u(t))\lesssim \mathcal{E}_2^\frac12+t\|\mathcal{L} u(t)\|_{L^2}.
\end{align*}
\end{lem}
The following lemma shows that $\|\langle c_a t-r\rangle P_a\partial_t^2\Gamma ^\alpha u\|_{L^2}$ can be controlled by the weighted norm $\mathcal{X}_{|\alpha|+2}$.
\begin{lem}\label{lem43}
  Suppose $ u\in H_\Gamma^k(T) $ is a solution of \eqref{1.2}, then there holds
  \begin{equation}\label{40}
 \sum_{|\alpha|\leq k-2}\|\langle c_a t-r\rangle P_a\partial_t^2\Gamma ^\alpha u\|_{L^2}\lesssim \mathcal{X}_k+\mathcal{X}_{[\frac{k-1}2]+3}\mathcal{E}_{k-1}^\frac12+\mathcal{X}_{k}\mathcal{E}_{[\frac{k-1}2]+3}^\frac12.
  \end{equation}
\end{lem}
{\bf Proof.} From \eqref{3.1}, we have
  \begin{equation}\label{41}
  \begin{split}
 \sum_{|\alpha|\leq k-2}\|\langle c_a t-r\rangle P_a\partial_t^2\Gamma ^\alpha u\|_{L^2}\leq&\sum_{|\alpha|\leq k-2}\|\langle c_a t-r\rangle P_a\nabla^2\Gamma ^\alpha u\|_{L^2}\\
 &+C\sum_{|\alpha|\leq k-2}\sum_{|\beta+\gamma|=|\alpha|}\|\langle c_a t-r\rangle P_a\nabla^2\Gamma ^\beta u\nabla\Gamma^\gamma u\|_{L^2}\\
 &+\sum_{|\alpha|\leq k-2}\sum_{|\beta+\gamma|=|\alpha|}\|\langle c_a t-r\rangle P_a\Lambda ^\beta \tilde{\rho}\partial_t^2\Gamma^\gamma u\|_{L^2}.
 \end{split}
  \end{equation}
 If $\|\tilde{\rho}\|_{H_\Lambda^k}\leq\frac12 $, then the last term of \eqref{41} can be absorbed by the left. Then we have
     \begin{equation}\label{42}
 \sum_{|\alpha|\leq k-2}\|\langle c_a t-r\rangle P_a\partial_t^2\Gamma ^\alpha u\|_{L^2}\lesssim \mathcal{X}_k+\sum_{|\alpha|\leq k-2}\sum_{|\beta+\gamma|=|\alpha|}\|\langle c_a t-r\rangle P_a\nabla^2\Gamma ^\beta u\nabla\Gamma^\gamma u\|_{L^2}
  \end{equation}
  Because that $|\beta+\gamma|\leq k-2$, either $|\beta|\leq [\frac{k-1}2]-1$ or $|\gamma|\leq [\frac{k-1}2]$ holds. If $|\beta|\leq [\frac{k-1}2]-1$, then by \eqref{4.3} in Lemma \ref{lem4.1},
\begin{equation}\label{43}
\begin{split}
&\sum_{|\alpha|\leq k-2}\sum_{|\beta+\gamma|=|\alpha|}\|\langle c_a t-r\rangle P_a\nabla^2\Gamma ^\beta u\nabla\Gamma^\gamma u\|_{L^2}\\
\leq& \sum_{|\alpha|\leq k-2}\sum_{|\beta+\gamma|=|\alpha|}\|\langle r\rangle\langle c_a t-r\rangle P_a\nabla^2\Gamma ^\beta u\|_{L^\infty}\|\nabla\Gamma^\gamma u\|_{L^2}\\
\lesssim &\mathcal{X}_{[\frac{k-1}2]+3}\mathcal{E}_{k-1}^\frac12.
 \end{split} \end{equation}
 If $|\gamma|\leq [\frac{k-1}2]$, then by \eqref{4.1} in Lemma \ref{lem4.1},
 \begin{equation}\label{44}
\begin{split}
&\sum_{|\alpha|\leq k-2}\sum_{|\beta+\gamma|\leq|\alpha|}\|\langle c_a t-r\rangle P_a\nabla^2\Gamma ^\beta u\nabla\Gamma^\gamma u\|_{L^2}\\
\leq& \sum_{|\alpha|\leq k-2}\sum_{|\beta+\gamma|\leq|\alpha|}\|\langle c_a t-r\rangle P_a\nabla^2\Gamma ^\beta u\|_{L^2}\|\nabla\Gamma^\gamma u\|_{L^\infty}\\
\lesssim &\mathcal{X}_{k}\mathcal{E}_{[\frac{k-1}2]+3}^\frac12.
 \end{split} \end{equation}
Finally, \eqref{40} is completed from \eqref{42}-\eqref{44}.
\hfill$\Box$
\begin{lem}\label{lem4.2}
Suppose $ u\in H_\Gamma^k(T) $ is a solution of the equation \eqref{1.2}, then
\begin{equation}\label{4.4}
  \mathcal{X}_k\lesssim\mathcal{E}_k^\frac12+ \mathcal{X}_{[\frac{k-1}2]+3} \mathcal{E}_{k-1}^\frac12+ \mathcal{X}_k \mathcal{E}_{[\frac{k-1}2]+3}^\frac12.
\end{equation}
\end{lem}
{\bf Proof.} Applying Lemma \ref{lem41}, we have that
\begin{align*}
  \mathcal{X}_k(u(t))\lesssim \mathcal{E}_k^\frac12+\sum_{|\alpha|\leq k-2}t\|\mathcal{L}\Gamma^\alpha u(t)\|_{L^2}.
\end{align*}
By \eqref{3.1}, we obtain that
\begin{equation}\label{4.6}
  \mathcal{X}_k(u(t))\lesssim \mathcal{E}_k^\frac12+\sum_{|\beta+\gamma|\leq k-2}t\|\nabla^2\Gamma^\beta u\nabla^\gamma u\|_{L^2}
  +\sum_{|\beta+\gamma|\leq k-2}t\|\Lambda^\beta\tilde{\rho} \partial_t^2\Gamma^\gamma u\|_{L^2}.
\end{equation}
For $t\|\nabla^2\Gamma^\beta u\nabla^\gamma u\|_{L^2}$, there holds either $|\beta|\leq [\frac{k-1}2]-1$ or $|\gamma|\leq [\frac{k-1}2]$. Note that
 \begin{align}\label{4.7}
  \langle t\rangle^{-1}\langle r\rangle\sum_a\langle c_a t-r\rangle P_a \gtrsim I,
 \end{align}
 thus we can control it by \begin{align*}
  t\|\nabla^2\Gamma^\beta u\nabla^\gamma u\|_{L^2}\lesssim\sum_a\|\langle r\rangle\langle c_a t-r\rangle P_a\nabla^2\Gamma^\beta u\|_{L^\infty}\|\nabla \Gamma^\gamma u\|_{L^2}
\end{align*}
while $|\beta|\leq [\frac{k-1}2]-1$, and
\begin{align*}
  t\|\nabla^2\Gamma^\beta u\nabla^\gamma u\|_{L^2}\lesssim\sum_a\|\langle c_a t-r\rangle P_a\nabla^2\Gamma^\beta u\|_{L^2}\|\langle r\rangle\nabla \Gamma^\gamma u\|_{L^\infty}
\end{align*}
while $|\gamma|\leq [\frac{k-1}2]$.\\
Because of  Lemma \ref{lem4.1}, we obtain
\begin{align}\label{4.8}
 \sum_{|\beta+\gamma|\leq k-2} t\|\nabla^2\Gamma^\beta u\nabla^\gamma u\|_{L^2}\lesssim  \mathcal{X}_{[\frac{k-1}2]+3} \mathcal{E}_{k-1}^\frac12+ \mathcal{X}_k \mathcal{E}_{[\frac{k-1}2]+3}^\frac12.
\end{align}
For the last term of \eqref{4.6}, note \eqref{4.7}, we have
\begin{equation*}
\begin{split}
 \sum_{|\beta+\gamma|\leq k-2}t\|\Lambda^\beta\tilde{\rho} \partial_t^2\Gamma^\gamma u\|_{L^2}&\lesssim \sum_a\sum_{|\beta+\gamma|\leq k-2}\|\langle r\rangle\langle c_a t-r\rangle P_a\Lambda^\beta\tilde{\rho}\partial_t^2\Gamma^\gamma u\|_{L^2}\\
 &\lesssim\sum_a\sum_{|\beta+\gamma|\leq k-2}\|\langle r\rangle\Lambda^\beta\tilde{\rho}\|_{L^\infty}\|\langle c_a t-r\rangle P_a\partial_t^2\Gamma^\gamma u\|_{L^2},
 \end{split}
\end{equation*}
then by \eqref{21} and Lemma \ref{lem43},
\begin{equation}\label{4.14}
\sum_{|\beta+\gamma|\leq k-2}t\|\Lambda^\beta\tilde{\rho} \partial_t^2\Gamma^\gamma u\|_{L^2}\lesssim\delta\big( \mathcal{X}_k+\mathcal{X}_{k}\mathcal{E}_{[\frac{k-1}2]+3}^\frac12+\mathcal{X}_{[\frac{k-1}2]+3}\mathcal{E}_{k-1}^\frac12\big).
\end{equation}
 Concluding \eqref{4.6}, \eqref{4.8} and \eqref{4.14}, there exists a uniform constant $C$ such that
\begin{equation}\label{4.15}
  \mathcal{X}_k(u(t))\leq C \Big(\mathcal{E}_k^\frac12+ \mathcal{X}_{[\frac{k-1}2]+3} \mathcal{E}_{k-1}^\frac12+ \mathcal{X}_k \mathcal{E}_{[\frac{k-1}2]+3}^\frac12+\delta( \mathcal{X}_k+\mathcal{X}_{k}\mathcal{E}_{[\frac{k-1}2]+3}^\frac12+\mathcal{X}_{[\frac{k-1}2]+3}\mathcal{E}_{k-1}^\frac12)\Big).
\end{equation}
If $\delta\leq\frac1{2C} $, then we have
\begin{equation*}
  \mathcal{X}_k(u(t))\leq (2C+1) \Big(\mathcal{E}_k^\frac12+ \mathcal{X}_{[\frac{k-1}2]+3} \mathcal{E}_{k-1}^\frac12+ \mathcal{X}_k \mathcal{E}_{[\frac{k-1}2]+3}^\frac12\Big).
\end{equation*}
\hfill$\Box$

\begin{lem}\label{lem4.3}
 Let $k\geq 9$. Suppose $ u\in H_\Gamma^k(T) $ is a solution of the equation\ \eqref{1.2}, if $\mathcal{E}_{k-2}$ is small enough for $0\leq t\leq T$, then
\begin{equation}\label{4.4}
\mathcal{X}_{k-2} (u(t))\lesssim{\mathcal{E}}_{k-2}^{\frac{1}{2}}(u(t)),
\end{equation}
\begin{equation}\label{4.5}
\mathcal{X}_k(u(t))\lesssim\mathcal{E}_k^{\frac{1}{2}}(u(t))[1+\mathcal{E}_{k-2} ^{\frac{1}{2}}(u(t))],
\end{equation}
and
\begin{equation}\label{4.9}
 \sum_{|\alpha|\leq k-2}\|\langle c_a t-r\rangle P_a\partial_t^2\Gamma ^\alpha u\|_{L^2}\lesssim\mathcal{X}_k(u(t))\lesssim\mathcal{E}_k^{\frac{1}{2}}(u(t)).
\end{equation}

\end{lem}
The proof of Lemma \ref{lem4.3} is evident by Lemma \ref{lem4.2} and Lemma \ref{lem43}.

\section{Almost global existence}
In this section, we will show the higher-order energy estimate
\begin{equation}\label{51}
\mathcal{E}'_{k}(u(t))\lesssim {\langle t \rangle}^{-1}\mathcal{E}_{k}(u(t))\mathcal{E}_{k-2}(u(t))^{\frac{1}{2}}+\langle t \rangle^{-1}\delta \mathcal{E}_{k}(u(t)),
\end{equation}
which can be adapted to prove the solution of \eqref{1.2} exists almost globally.

Suppose that $u(t)\in H_\Gamma^k(T)$ is a local solution of \eqref{1.2}. Taking inner product with $\partial_t \Gamma ^\alpha u$ on both sides of \eqref{3.1}, we get
\begin{equation}\label{52}
\begin{split}
&\frac{1}{2} \frac{d}{dt}\int_{\mathbb{R}^3}|\partial_t\Gamma^\alpha u|^2+c_2^2|\triangledown  \Gamma^\alpha u|^2 + (c_1^2-c_2^2)(\triangledown\cdot\Gamma^\alpha u)^2dx\\
\backsimeq&\sum_{|\beta+\gamma|=|\alpha|}\int_{\mathbb{R}^3}(\partial_t\Gamma^\alpha u,N(\Gamma^\beta u,\Gamma^\gamma u))dx -\sum_{|\beta+\gamma|= |\alpha|}\int_{\mathbb{R}^3}\Lambda^\beta\tilde{\rho} (\partial_t\Gamma^\alpha u,\partial_t^2\Gamma^\gamma u) dx.
\end{split}
\end{equation}
For the first term on the right side of \eqref{52}, if $\beta=\alpha $ or $\gamma= \alpha$, we have
\begin{equation*}
\begin{split}
\int_{\mathbb{R}^3}(\partial_t\Gamma^\alpha u,N(\Gamma^\alpha u,u))dx=&B_{lmn}^{ijk}\int_{\mathbb{R}^3}\partial_t\Gamma^\alpha u^i \partial_l(\partial_m\Gamma^\alpha u^j\partial_n u^k)dx\\
=&-B_{lmn}^{ijk}\int_{\mathbb{R}^3}\partial_t\partial_l\Gamma^\alpha u^i \partial_m\Gamma^\alpha u^j\partial_n u^k dx\\
=&-\frac{1}{2}\frac{d}{dt}B_{lmn}^{ijk}\int_{\mathbb{R}^3}\partial_l\Gamma^\alpha u^i \partial_m \Gamma^\alpha u^j \partial_n u^k dx\\
&+\frac{1}{2}B_{lmn}^{ijk}\int_{\mathbb{R}^3}\partial_l\Gamma^\alpha u^i\partial_m\Gamma^\alpha u^j\partial_t\partial_n u^k dx.
\end{split}
\end{equation*}
We set
$$\tilde{N}(u,v,w)=B_{lmn}^{ijk}\partial_l u^i \partial_m v^j \partial_n w^k ,$$
and
\begin{equation}\label{5.1}
\tilde{\mathcal{E}_{k}}(u(t))=\mathcal{E}_{k}(u(t))+\sum\limits_{|\alpha|\leq k-1}\int_{\mathbb{R}^3}\tilde{N}(\Gamma^\alpha u,\Gamma^\alpha u,u)dx.
\end{equation}
Thus
\begin{equation}\label{5.2}
\begin{split}
\tilde{\mathcal{E}}_{k}'(u(t))\backsimeq&\sum_{|\alpha|\leq k-1}\sum\limits_{|\beta|+|\gamma|= |\alpha|\atop |\beta|,|\gamma|\neq |\alpha|}\|\partial \Gamma^\alpha u\|_{L^2}\|\partial\nabla\Gamma^\beta u\nabla\Gamma^\gamma u\|_{L^2}\\
&+\sum_{|\alpha|\leq k-1}\sum_{|\beta+\gamma|= |\alpha|}\int_{\mathbb{R}^3}\Lambda^\beta\tilde{\rho} \partial_t\Gamma^\alpha u\partial_t^2\Gamma^\gamma u dx.
\end{split}
\end{equation}
Since $|\beta|+|\gamma|\leq |\alpha|$ and $|\beta|,|\gamma|\neq |\alpha|$, we may assume either $|\beta|+1\leq [\frac k2]$ or $|\gamma|\leq [\frac k2]$. As the way of proving Lemma \ref{lem4.2} and by Lemma \ref{lem4.1} and Lemma \ref{lem4.3}, we have that
\begin{equation}\label{5.5}
  \sum_{|\alpha|\leq k-1,}\sum\limits_{|\beta|+|\gamma|= |\alpha|\atop |\beta|,|\gamma|\neq |\alpha|}\|\partial \Gamma^\alpha u\|_{L^2}\|\partial\nabla\Gamma^\beta u\nabla\Gamma^\gamma u\|_{L^2}\lesssim \langle t \rangle^{-1}\mathcal{E}_{k-2}^{\frac{1}{2}}\mathcal{E}_k.
\end{equation}
To estimate the second term on the right side of \eqref{5.2}, we revisit it as
\begin{equation}\label{5.6}
\begin{split}
\sum_{|\alpha|\leq k-1}\sum_{|\beta+\gamma|= |\alpha|}\int_{\mathbb{R}^3}\Lambda^\beta\tilde{\rho} \partial_t\Gamma^\alpha u\partial_t^2\Gamma^\gamma u dx=&\sum_{|\alpha|\leq k-1}\int_{\mathbb{R}^3}\tilde{\rho} \partial_t\Gamma^\alpha u\partial_t^2\Gamma^\alpha u dx\\
&+\sum_{|\alpha|\leq k-1}\sum_{|\beta+\gamma|= |\alpha|\atop |\gamma|\leq k-2}\int_{\mathbb{R}^3}\Lambda^\beta\tilde{\rho} \partial_t\Gamma^\alpha u\partial_t^2\Gamma^\gamma u dx.
\end{split}
\end{equation}
Since
\begin{equation*}
 \sum_{|\alpha|\leq k-1}\int_{\mathbb{R}^3}\tilde{\rho} \partial_t\Gamma^\alpha u\partial_t^2\Gamma^\alpha u dx=- \sum_{|\alpha|\leq k-1}\frac12\frac d{dt}\int_{\mathbb{R}^3}\tilde{\rho} |\partial_t\Gamma^\alpha u|^2 dx,
\end{equation*}
we set
  \begin{equation}\label{5.8}
\hat{\mathcal{E}_{k}}(u(t))=\tilde{\mathcal{E}_{k}}(u(t))+\sum_{|\alpha|\leq k-1}\int_{\mathbb{R}^3}\tilde{\rho} |\partial_t\Gamma^\alpha u|^2 dx,
\end{equation}
then we have that
\begin{equation}\label{5.9}
\begin{split}
\hat{\mathcal{E}_{k}}'(u(t)) \lesssim \langle t \rangle^{-1}\mathcal{E}_{k-2}^{\frac{1}{2}}\mathcal{E}_k+\sum_{|\alpha|\leq k-1}\sum_{|\beta+\gamma|= |\alpha|\atop |\gamma|\leq k-2}\int_{\mathbb{R}^3}\big|\Lambda^\beta\tilde{\rho} \partial_t\Gamma^\alpha u\partial_t^2\Gamma^\gamma u \big|dx
  \end{split}
  \end{equation}
combined with \eqref{5.2}-\eqref{5.8}. By \eqref{21}, \eqref{4.7} and \eqref{4.9}, we have
\begin{equation}
\begin{split}
&\sum_{|\alpha|\leq k-1}\sum_{|\beta+\gamma|= |\alpha|\atop |\gamma|\leq k-2}\int_{\mathbb{R}^3}\big|\Gamma^\beta\tilde{\rho} \partial_t\Gamma^\alpha u\partial_t^2\Gamma^\gamma u \big|dx\\
\lesssim&\sum_a\sum_{|\alpha|\leq k-1}\sum_{|\gamma|\leq k-2} \langle t\rangle^{-1}\delta\|\langle c_at-r\rangle P_a\partial_t^2\Gamma^\gamma u \|_{L^2}\|\partial_t\Gamma^\alpha u\|_{L^2}\\
\lesssim&\langle t \rangle^{-1}\delta \mathcal{X}_{k}\mathcal{E}_{k}^{\frac12}.
\end{split}
\end{equation}
Consequently, one has
\begin{equation}
  \hat{\mathcal{E}_{k}}'(u(t)) \lesssim \langle t \rangle^{-1}\mathcal{E}_{k-2}^{\frac{1}{2}}\mathcal{E}_k+\langle t \rangle^{-1}\delta \mathcal{E}_{k},
\end{equation}
which gives that
\begin{equation}\label{511}
  \mathcal{E}_{k}(u(t))\lesssim M^2\langle t \rangle^{\varepsilon+\delta},
\end{equation}
because the modified energy $\hat{\mathcal{E}_{k}}(u(t))$ is equivalent to the standard one $\mathcal{E}_{k}(u(t))$ and provided that $\mathcal{E}_{k-2}(u(t))\lesssim\mathcal{E}_{k-2}(u(0))$, which implies the almost global existence result of John \cite{JF} (see also \cite{KS}, \cite{S2}).

\section{Global existence}
This section is devoted to prove the existence of the global solution of the \eqref{1.2} under assumptions of Theorem \ref{thm 2.2}. Combined with the higher-order energy estimate, we will discuss the lower-order energy estimate
\begin{equation}\label{6.1}
\mathcal{E}'_{k-2}(u(t))\lesssim{\langle t \rangle}^{-\frac32}\mathcal{E}_{k}^{\frac{1}{2}}(u(t))\mathcal{E}_{k-2}(u(t))+ {\langle t \rangle}^{-\frac32}\delta\mathcal{E}_{k}^{\frac{1}{2}}(u(t))\mathcal{E}_{k-2}^{\frac{1}{2}}(u(t)).
\end{equation}

From \eqref{52}, \eqref{5.1} and \eqref{5.8}, let $|\alpha|\leq k-3$ , we have
\begin{equation}\label{6.2}
\begin{split}
\hat{\mathcal{E}}_{k-2}'(u(t))\backsimeq&\sum_{|\alpha|\leq k-3}\sum\limits_{|\beta|+|\gamma|= |\alpha|\atop |\beta|,|\gamma|\neq |\alpha|}\|\partial \Gamma^\alpha u\|_{L^2}\|\partial\nabla\Gamma^\beta u\nabla\Gamma^\gamma u\|_{L^2}\\
&+\sum_{|\alpha|\leq k-3}\sum_{|\beta+\gamma|= |\alpha|\atop |\gamma|\leq k-4}\int_{\mathbb{R}^3}\Lambda^\beta\tilde{\rho} \partial_t\Gamma^\alpha u\partial_t^2\Gamma^\gamma u dx.
\end{split}
\end{equation}
Like the discussion of Sideris,  away from the inner light cone, the first term on the right side of \eqref{6.2} can be controlled by $\langle t \rangle ^{-\frac32}\mathcal{E}_{k-2}\mathcal{E}_k^\frac12$ directly by Lemma \ref{lem4.1}.  When $r$ is comparable to $\langle t\rangle$,  it also can be  achieved successfully thanks to the null condition \eqref{1.5} and \eqref{1.6}. We omit the details here, readers can refer to \cite{S2} to get
\begin{align}\label{6.3}
\sum_{|\alpha|\leq k-3}\sum\limits_{|\beta|+|\gamma|= |\alpha|\atop |\beta|,|\gamma|\neq |\alpha|}\|\partial \Gamma^\alpha u\|_{L^2}\|\partial\nabla\Gamma^\beta u\nabla\Gamma^\gamma u\|_{L^2} \lesssim \langle t \rangle ^{-\frac{3}{2}}\mathcal{E}_k^{\frac{1}{2}}\mathcal{E}_{k-2}.
\end{align}
We focus on the last term of \eqref{6.2} in two cases: $\langle c_2 t\rangle\leq r$ and $\langle c_2 t\rangle\geq r$.  Note that \eqref{21}, \eqref{4.7}, \eqref{4.9} and \eqref{4.2-1}, we have
\begin{equation}\label{6.3}
\begin{split}
&\sum_{|\alpha|\leq k-3}\sum_{|\beta+\gamma|= |\alpha|\atop |\gamma|\leq k-4}\int_{r\geq\langle c_2 t\rangle}\Lambda^\beta\tilde{\rho} \partial_t\Gamma^\alpha u\partial_t^2\Gamma^\gamma u dx\\
\lesssim& \sum_a\sum_{|\alpha|\leq k-3}\sum_{|\beta+\gamma|= |\alpha|\atop |\gamma|\leq k-4}\langle t\rangle^{-1}\|\langle r\rangle\Lambda^\beta\tilde{\rho}\|_{L^2}\| \partial_t\Gamma^\alpha u\|_{L^\infty( r\geq\langle c_2 t\rangle)}\|\langle c_at-r\rangle P_a\partial_t^2\Gamma^\gamma u\|_{L^2}\\
\lesssim&\sum_{|\alpha|\leq k-3}\langle t\rangle^{-1}\delta\| \partial_t\Gamma^\alpha u\|_{L^\infty( r\geq\langle c_2 t\rangle)}\mathcal{E}_{k-2}^\frac12\\
\lesssim&\sum_{|\alpha|\leq k-3}\langle t\rangle^{-2}\delta\mathcal{E}_{k}^\frac12\mathcal{E}_{k-2}^\frac12.
\end{split}
\end{equation}
When $\langle c_2 t\rangle\geq r$, the ratio $\langle t\rangle^{-1}\langle c_at-r\rangle$ and $\langle t\rangle^{-\frac12}\langle c_1t-r\rangle^\frac12$ are bounded below, with \eqref{4.2} in Lemma \ref{lem4.1} thus
\begin{equation}\label{6.4}
\begin{split}
&\sum_{|\alpha|\leq k-3}\sum_{|\beta+\gamma|= |\alpha|\atop |\gamma|\leq k-4}\int_{r\leq\langle c_2 t\rangle}\Lambda^\beta\tilde{\rho} \partial_t\Gamma^\alpha u\partial_t^2\Gamma^\gamma u dx\\
\lesssim& \sum_a\sum_{|\alpha|\leq k-3}\sum_{|\beta+\gamma|= |\alpha|\atop |\gamma|\leq k-4}\langle t\rangle^{-1}\|\Lambda^\beta\tilde{\rho}\|_{L^2}\| \partial_t\Gamma^\alpha u\|_{L^\infty( r\leq\langle c_2 t\rangle)}\|\langle c_at-r\rangle P_a\partial_t^2\Gamma^\gamma u\|_{L^2}\\
\lesssim&\sum_{|\alpha|\leq k-3}\langle t\rangle^{-\frac32}\delta\|\langle c_1t-r\rangle^\frac12 \partial_t\Gamma^\alpha u\|_{L^\infty}\mathcal{E}_{k-2}^\frac12\\
\lesssim&\sum_{|\alpha|\leq k-3}\langle t\rangle^{-\frac32}\delta\mathcal{E}_{k}^\frac12\mathcal{E}_{k-2}^\frac12,
\end{split}
\end{equation}
Consequently, one has
 \begin{equation*}
\mathcal{E}'_{k-2}(u(t))\lesssim{\langle t \rangle}^{-\frac32}\mathcal{E}_{k}^{\frac{1}{2}}(u(t))\mathcal{E}_{k-2}(u(t))+ {\langle t \rangle}^{-\frac32}\delta\mathcal{E}_{k}^{\frac{1}{2}}(u(t))\mathcal{E}_{k-2}^{\frac{1}{2}}(u(t)).
\end{equation*}
Attached with \eqref{511}, we have
 \begin{equation*}
\frac{d\mathcal{E}_{k-2}^\frac12(u(t))}{dt}\lesssim\frac12{\langle t \rangle}^{\frac{-3+\varepsilon+\delta}2}M\mathcal{E}_{k-2}^\frac12(u(t))+ \frac12\delta M{\langle t \rangle}^{\frac{-3+\varepsilon+\delta}2}.
\end{equation*}
Multiplying $e^{-\frac M{-1+\varepsilon+\delta}\langle t\rangle^{\frac{-1+\varepsilon+\delta}2}}$ on both sides of the above, and solve the ordinary differential inequality, there holds
\begin{align*}
  \mathcal{E}_{k-2}^\frac12(u(t))&\lesssim e^{\frac M{-1+\varepsilon+\delta}\langle t\rangle^{\frac{-1+\varepsilon+\delta}2}}\varepsilon e^{C_2M}+\delta M\int_0^t e^{-\frac M{-1+\varepsilon+\delta}\langle s\rangle^{\frac{-1+\varepsilon+\delta}2}}\langle s\rangle^{\frac{-3+\varepsilon+\delta}2}ds\\
  &\lesssim\varepsilon e^{C_2M}+\delta e^{C_3M},
\end{align*}
provided $\varepsilon+\delta<<1$, where $C_2$ and $C_3$ are two constants uniformly in $t$. The proof of Theorem \ref{thm 2.2} is completed.

\section*{Acknowledgemets}
The first author was supported by China Postdoctoral Science Foundation funded project (D.10-0101-17-B02).

\end{document}